\documentclass[12pt,sort&compress]{elsarticle}
\usepackage{amsmath}
\usepackage{amsfonts}
\usepackage{amssymb}
\usepackage{amsthm}
\usepackage{bm}
\usepackage{indentfirst}
\usepackage{graphicx}
\usepackage{subfigure}
\usepackage{array}
\usepackage{fullpage}
\usepackage{color}

\DeclareMathAlphabet{\mathsfsl}{OT1}{cmss}{m}{sl}

\newcommand{\PreserveBackslash}[1]{\let\temp=\\#1\let\\=\temp}
\newcolumntype{C}[1]{>{\PreserveBackslash\centering}p{#1}}
\newcolumntype{R}[1]{>{\PreserveBackslash\raggedleft}p{#1}}
\newcolumntype{L}[1]{>{\PreserveBackslash\raggedright}p{#1}}

\renewcommand{\vartheta}{\Theta}

\definecolor{mygreen}{rgb}{0.1,0.75,0.2}

\numberwithin{equation}{section}

\theoremstyle{definition}

\usepackage{hyperref}
\usepackage{fancyvrb}
\usepackage{fancyhdr}
\biboptions{sort&compress}

\usepackage{bm}
\usepackage{mathrsfs}
\usepackage{amsmath}
\usepackage{amsfonts}
\usepackage{amsthm}
\usepackage{color}
\usepackage{setspace}
\usepackage{float}
\usepackage{url}
\usepackage{multicol}
\usepackage{subfigure}
\usepackage{bm}
\usepackage[top=1in,bottom=1in,left=1in,right=1in]{geometry}

\begin{document}

\theoremstyle{definition}
\newtheorem{remark}{Remark}

\newcommand{\vertiii}[1]{{\left\vert\left\vert\left\vert #1
    \right\vert\right\vert\right\vert}}

\begin{frontmatter}

\title{Asymptotically compatible meshfree discretization of state-based peridynamics for linearly elastic composite materials}


\author[nt]{Nathaniel Trask \footnote{Sandia National Laboratories is a multimission laboratory managed and operated by National Technology and Engineering Solutions of Sandia, LLC.,a wholly owned subsidiary of Honeywell International, Inc., for the U.S. Department of Energy’s National Nuclear Security Administration under contract DE-NA-0003525. This paper describes objective technical results and analysis. Any subjective views or opinions that might be expressed in the paper do not necessarily represent the views of the U.S. Department of Energy or the United States Government.}}\ead{natrask@sandia.gov}
\cortext[cor1]{Corresponding Author}
\author[bh]{Benjamin Huntington}
\author[dl]{David Littlewood}

%
\begin{abstract}
State-based peridynamic models provide an important extension of bond-based models that allow the description of general linearly elastic materials. Meshfree discretizations of these nonlocal models are attractive due to their ability to naturally handle fracture. However, singularities in the integral operators have historically proven problematic when seeking convergent discretizations. We utilize a recently introduced optimization-based quadrature framework to obtain an asymptotically compatible scheme able to discretely recover local linear elasticity as the nonlocal interaction is reduced at the same rate as the grid spacing. By introducing a correction to the definition of nonlocal dilitation, surface effects for problems involving bond-breaking and free surfaces are avoided without the need to modify the material model. We use a series of analytic benchmarks to validate the consistency of this approach, illustrating second-order convergence for the Dirichlet problem, and first-order convergence for problems involving curvilinear free surfaces and composite materials. We additionally illustrate that these results hold for material parameters chosen in the near-incompressible limit.
\end{abstract}

\begin{keyword}
Nonlocal \sep Peridynamics \sep asymptotic compatibility \sep meshfree \sep incompressibility 

\end{keyword}

\end{frontmatter}

\section{Introduction}

Consider a compact connected domain $\Omega$ with Lipschitz continuous boundary $\partial \Omega$. We seek a solution on this domain to the linear elasticity equations. 
\begin{align}\label{eq:locElasticity}
  \begin{cases}
    -\nabla \cdot \sigma = \mathbf{f}\\
    \sigma = \lambda \left(\nabla \cdot \mathbf{u} \right) \mathbf{I} + \mu \left(\nabla \mathbf{u} + \nabla \mathbf{u}^\intercal \right)\\
    \mathbf{u}|_{\partial\Omega_D} = \mathbf{u}_D\\
    \hat{\mathbf{n}} \cdot \sigma |_{\partial\Omega_N} = 0\\
  \end{cases}
\end{align}
where $\lambda$ and $\mu \in L_2(\Omega)$ are Lam\`e coefficients, and the boundary is partitioned into Dirichlet and Neumann parts ($\partial \Omega = \partial \Omega_D \cup \partial \Omega_N$), in which $\mathbf{u}_D$ denotes a prescribed displacement condition and a traction-free condition is imposed on the Neumann part. In this work we specifically consider problems involving fracture of composite materials constituted of $n$ phases, so that the domain may be partitioned into a collection of disjoint subdomains with piecewise constant material properties in each subdomain, i.e. $\Omega = \underset{i}{\cup}{\Omega_i}$, $\Omega_i\cap \Omega_j = \emptyset$, and $\lambda = \lambda_i$, $\mu = \mu_i$ for $\mathbf{x} \in \Omega_i$.

As fracture and discontinuous material properties both lead naturally to solutions with reduced regularity, we will discretize this problem with a nonlocal theory of mechanics which poses the governing equations in terms of compactly supported integral operators, so that the problem remains well-posed in the presence of discontinuities. Specifically, we consider here the state-based linear peridynamic solid model \cite{silling2010linearized}. This model consists of a nonlocal parameter $\delta$ characterizing the degree of nonlocality, and in the limit as $\delta \rightarrow 0$ the model recovers Equation \ref{eq:locElasticity} for general choice of $\mu$ and $\lambda$. In this work we thus consider nonlocal models used as a \textit{regularization} of classical continuum mechanics, and we seek solutions to this nonlocal model which discretely converge to a physically relevant local counterpart. In contrast to classical local models, cracks emerge naturally from the solution\cite{silling2010crack}, avoiding the need for computationally involved techniques popular in the finite element community, such as tracking of free-surfaces and enrichment strategies.

For a given quasi-uniform discretization of a nonlocal model, there exists a representative lengthscale $h$ characterizing the discretization resolution. To maintain an easily scalable implementation, it is desirable that the ratio $\frac{h}{\delta}$ be bound by a constant $M$ as $\delta \rightarrow 0$. This so-called "M-convergence" results in a sparse linear system with bounded bandwidth that may be solved efficiently with standard preconditioning techniques \cite{bobaru2009convergence}. Convergence of the nonlocal solution to the local counterpart under this M-convergence condition is particularly challenging for nonlocal models due to the presence of singular integrands\cite{trask2019asymptotically}, surface effects near boundaries\cite{le2017surface,mitchell2013towards}, and issues with handling internal interfaces\cite{seleson2013interface}; methods which are able to achieve this convergence are termed asymptotically compatible (AC)\cite{tian_2014}. We therefore restrict the focus of this paper to nonlocal models as a computationally efficient regularization of local models, as opposed to the discretization of models describing nonlocal physical processes, where $\delta$ represents a physically nonlocal interaction and may not be refined with $h$ (see e.g., \cite{eringen1972nonlocal,silling_2000}).

In contrast to variational discretizations of nonlocal operators, meshfree discretizations are often used due to their ease of implementation, efficiency, and natural treatment of fracture problems \cite{silling_2005_2}. In a previous work, we have shown that a meshfree optimization-based quadrature (OBQ) rule may be used to efficiently guarantee asymptotic compatibility for nonlocal operators of general form\cite{trask2019asymptotically}. In this work, the prototypical microelastic brittle (PMB) model \cite{silling_2005_2} was considered, which is an example of a bond-based material model. Unlike state-based nonlocal models, bond-based models are restricted to linearly elastic materials of Poisson ratio $\nu = \frac14$ in 3D, and thus to consider general materials the previous approach must be extended. We will show that while under certain conditions the previous approach works without modification, additional considerations must be taken to handle both bond-breaking damage models in a state-based setting, and also to stably handle jumps in material properties. We will show that handling these challenges will also remedy the so-called "surface effect" issues arising in the literature\cite{le2017surface,mitchell2013towards}. In our approach, we modify the definition of dilitation to ensure a consistent limit to the local divergence; this is in contrast to previous works which focus on adding stiffness in the vicinity of the boundary \cite{le2017surface}. Thus, we illustrate that an accurate treatment of quadrature and consistent definition of dilitation is sufficient to obtain an AC discretization of peridynamics without modification to the material model or physics.

We organize this paper as follows. We will define our scheme in Section \ref{sec:model}, after recalling first the model definition in Subsection \ref{ssec:LPS} and the meshfree quadrature scheme in Subsection \ref{ssec:OBQ}. We then introduce a new approach to preserve asymptotic compatibility as bonds are broken to model static loading of cracks in Subsection \ref{ssec:bondCorrection}. We proceed to investigate a number of two-dimensional statics problems with analytic solutions in Section \ref{sec:numerics}. We categorize these results as: manufactured solutions to illustrate convergence (Subsection \ref{ssec:manSol}); analytic solutions of homogeneous materials with free surfaces and cracks (Subsection \ref{ssec:isoSol}); and analytic solutions of composite materials with internal interfaces (Subsection \ref{ssec:anisoSol}).

\section{Model definition}\label{sec:model}

We recall both the LPS model and the OBQ technique for discretizing compactly supported nonlocal operators. For problems without fracture, we will see that this is sufficient to define an asymptotically compatible scheme. We also introduce a novel modification of the nonlocal dilitation that will preserve consistency for cracks modeled through bond-breaking techniques.

\subsection{Linear peridynamic solid model}\label{ssec:LPS}

We define the state-based linearized peridynamic material in $d$ dimensions as follows. Let $\theta$ be the nonlocal dilitation, generalizing the local divergence of displacement, and let $K(r)$ denote a positive radial function with compact support $\delta$. The momentum balance and nonlocal dilitation are then given by the following, for scaling parameters $C_\alpha$ and $C_\beta$.

\begin{equation}\label{eq:nonlocElasticity}
    \frac{C_\alpha}{m(\delta)}  \int_{B_\delta (\mathbf{x})} \left(\lambda(x,y) - \mu(x,y)\right) K(\left|\mathbf{y}-\mathbf{x}\right|) \left(\mathbf{y}-\mathbf{x} \right)\left(\theta(\mathbf{x}) + \theta(\mathbf{y}) \right) d\mathbf{y}
\end{equation}
\begin{equation*}
  +  \frac{C_\beta}{m(\delta)}\int_{B_\delta (\mathbf{x})} \mu(x,y) K(\left|\mathbf{y}-\mathbf{x}\right|)\frac{\left(\mathbf{y}-\mathbf{x}\right)\otimes\left(\mathbf{y}-\mathbf{x}\right)}{\left|\mathbf{y}-\mathbf{x}\right|^2} \cdot \left(\mathbf{u}(\mathbf{x}) - \mathbf{u}(\mathbf{y}) \right) d\mathbf{y} = \mathbf{f}
\end{equation*}
\begin{equation}\label{eq:nonlocDilitation}
  \theta(\mathbf{x}) =  \frac{d}{m(\delta)} \int_{B_\delta (\mathbf{x})} K(\left|\mathbf{y}-\mathbf{x}\right|) \left(\mathbf{y}-\mathbf{x}\right) \cdot \left(\mathbf{u}(\mathbf{y}) - \mathbf{u}(\mathbf{x}) \right) d\mathbf{y}
\end{equation}
where we define the weighted volume
\begin{equation}
m(\delta) = \int_{B_\delta (\mathbf{x})} K(\left|\mathbf{y}-\mathbf{x}\right|) |\mathbf{y}-\mathbf{x}|^2 d\mathbf{y}
\end{equation}
and adapt the notation $\mu(x,y)$ and $\lambda(x,y)$ to denote an average material property $p$ satisfying $p(z,z) = p(z)$. We will consider for the purposes of this work the harmonic mean 
\begin{equation}
   \frac{2}{ \mu(x,y)} =\frac{1}{\mu(x)} + \frac{1}{\mu(y)}.
\end{equation}

For $B_\delta(\mathbf{x})\cap \partial \Omega = \emptyset$ and $\mathbf{u} \in C^2(\Omega)$, this system may be shown to converge to Equation \ref{eq:locElasticity} as $\delta \rightarrow 0$ with appropriate choice of scaling parameters and weighting function $K$. While this limit holds for general choice of $K$, for concreteness we will in the current work restrict ourselves to
\begin{align}
K(r) =
\begin{cases}
\frac{1}{r} &\quad r < \delta\\
0           &\quad \text{else}.
\end{cases}
\end{align}

To recover parameters for 3D linear elasticity, we may select $C_\alpha = 3$, $C_\beta = 30$, and $m(\delta) = \pi \delta^4$. For 2D plane strain, we may select $C_\alpha = 2$, $C_\beta = 16$, and $m(\delta) = \frac{2 \pi \delta^3}{3}$. Confirmation that these scalings recover the local limit may be obtained via a straightforward yet tedious Taylor series calculation; Mathematica scripts may be found to complete this calculation in the online version of this paper.

In addition to limiting to the linear operators, it may easily be seen that the definition of $\theta(\mathbf{x})$ limits to a local divergence operator as $\delta \rightarrow 0$ when $\mathbf{u} \in C^1(\Omega)$ and $B_\delta(\mathbf{x})\cap \partial \Omega = \emptyset$. We sketch the proof in 2D and again provide a supporting calculation online. Under these assumptions the displacement admits the Taylor series expansion $\mathbf{u}(\mathbf{y}) = \mathbf{u}(\mathbf{x}) + \nabla \mathbf{u}(\mathbf{x}) \cdot \left( \mathbf{y} - \mathbf{x} \right) + O(\delta^2)$. By switching to polar coordinates, it follows from a symmetry argument that the only non-vanishing terms are $\partial_x u$ and $\partial_y v$, and after applying the scaling $m(\delta)$ one recovers the local divergence. Unfortunately, if the domain of integration is non-spherical due to proximity to the boundary, this symmetry argument no longer holds and an inconsistent dilitation is obtained. Thus, surface-effects manifest themselves in the definition of dilitation before any modeling assumptions are made regarding the material response. We will address this issue in Subsection \ref{ssec:bondCorrection} by introducing a new definition of dilitation guaranteed to discretely preserve an AC limit to the local divergence operator.

\subsection{Optimization-based meshfree quadrature}\label{ssec:OBQ}

Strong-form particle discretizations of peridynamics, favored for their computational simplicity and ease of handling fracture, pursue a discretization of the above system through the following one point quadrature rule at a collection of collocation points $X_h = \left\{\mathbf{x}_i\right\}_{i=1 \dots N}$\cite{silling_2010}.

\begin{equation}\label{eq:discreteNonlocElasticity}
    \frac{C_\alpha}{m(\delta)} \sum_{j \in B_\delta (\mathbf{x}_i)} \left(\lambda_{ij} - \mu_{ij}\right) K_{ij} \left(\mathbf{x}_j-\mathbf{x}_i\right)\left(\theta_i + \theta_j \right) \omega_j
\end{equation}
\begin{equation*}
  +  \frac{C_\beta}{m(\delta)}\sum_{j \in B_\delta (\mathbf{x}_i)} \mu_{ij} K_{ij}\frac{\left(\mathbf{x}_j-\mathbf{x}_i\right)\otimes\left(\mathbf{x}_j-\mathbf{x}_i\right)}{\left|\mathbf{x}_j-\mathbf{x}_i\right|^2} \cdot \left(\mathbf{u}_i- \mathbf{u}_j \right) \omega_j = \mathbf{f}_i
\end{equation*}
\begin{equation}\label{eq:discreteNonlocDilitation}
  \theta(\mathbf{x}) =  \frac{d}{m(\delta)} \sum_{j \in B_\delta (\mathbf{x}_i)}  K_{ij} \left(\mathbf{x}_j-\mathbf{x}_i\right) \cdot  \left(\mathbf{u}_j - \mathbf{u}_i \right) \omega_j
\end{equation}
where we adapt the notations $f(x_i) = f_i$ and $f(x_i,x_j)=f_{ij}$, and we specify $\left\{\omega_i\right\}$ as a to-be-determined collection of quadrature weights admitting interpretation as a measure associated with each collocation point. We define in this section an optimization-based approach to defining these weights extending previous work \cite{trask2019asymptotically}, constructed to ensure consistency guarantees.

We first characterize the distribution of collocation points as follows. Recall the definitions \cite{wendland2004scattered} of \textit{fill distance}
\begin{equation}
h_{{X}_h,\Omega} = \underset{\bm{x} \in \Omega}{\sup}\, \underset{x_i \in {X}_h}{\min}||x_i - x_j||_2,
\end{equation}
and \textit{separation distance}
\begin{equation}
q_{{X}_h} = \frac12 \underset{i \neq j}{\min} ||x_i - x_j|_.
\end{equation}
We assume that $X_h$ is \textit{quasi-uniform}, namely that there exists $c_{qu} > 0$ such that
\begin{equation}
q_{{X}_h} \leq h_{{X}_h,\Omega} \leq c_{qu} q_{{X}_h}.
\end{equation}

To summarize the approach outlined in Trask et al\cite{trask2019asymptotically}, we seek quadrature weights for integrals supported on balls of the form
\begin{equation}
I[f] := \int_{B_\delta (\mathbf{x_i})} f(x,y) dy \approx I_h[f] := \sum_{j \in B_\delta (\mathbf{x}_i)} f(x_i,x_j) \omega_{j,i}
\end{equation}
where we include the subscript $i$ in $\left\{\omega_{j,i}\right\}$ to denote that we seek a different family of quadrature weights for different subdomains $B_\delta(\mathbf{x}_i)$. We obtain these weights from the following optimization problem

\begin{align}\label{eq:quadQP}
  \underset{\left\{\omega_{j,i}\right\}}{\text{argmin}} \sum_{j \in B_\delta (\mathbf{x}_i)} \omega_{j,i}^2 \quad
  \text{such that}, \quad
  I_h[p] = I[p] \quad \forall p \in \mathbf{V}_h
\end{align}
where $\mathbf{V}_h$ denotes a Banach space of functions which should be integrated exactly. We refer to previous work\cite{trask2019asymptotically} for further information, analysis, and implementation details. Provided the quadrature points are unisolvent over the desired reproducing space, this problem may be proven to have a solution by interpreting it as a generalized moving least squares (GMLS) problem; a rigorous discussion of this fact requires some mathematical machinery and will be discussed in a future work. For certain choices of $V_h$, such as $m^{th}$-order polynomials, Problem \ref{eq:quadQP} may be proven to have a solution for domains $\Omega$ satisfying a cone condition, quasi-uniform pointsets, and a sufficiently large number of quadrature points\cite{wendland2004scattered}. In previous work we have provided error estimates establishing the dependence of the convergence rate upon the resulting quadrature rule as a function of the order singularity appearing in the integrand.

In this work, we choose a reproducing space sufficient to exactly integrate Equations \ref{eq:discreteNonlocElasticity} and \ref{eq:discreteNonlocDilitation} in the case where $\mathbf{u}$ and $\theta$ are vector-valued polynomials of order two and scalar-valued polynomials of order one, respectively. To do this, we select as reproducing space $V_h = P_2 \cup S_{K,2,\mathbf{x}}$, where $P_2$ denotes scalar-valued quadratic polynomials, and 
$$S_{K,2,\mathbf{x}} := \left\{ \frac{(\mathbf{y} - \mathbf{x})\otimes(\mathbf{y} - \mathbf{x})}{|\mathbf{y} - \mathbf{x}|^3} \cdot \left(\mathbf{p}(\mathbf{y}) - \mathbf{p}(\mathbf{x})\right) \, \mathbf{p} \in (P_2)^d \right\}$$
denotes the space of vector-valued quadratic polynomials integrated against the tensorial kernel appearing in Equation \ref{eq:nonlocElasticity}. 

As discussed in \cite{trask2019asymptotically}, the key to obtaining these quadrature weights is that they may be evaluated analytically. In the online version of the paper we include Mathematica scripts calculating analytical expressions for the equality constraints in Equation \ref{eq:quadQP}. Obtaining the desired quadrature weights from the above equality constrained quadratic program thus amounts to solving a small linear system at each collocation point. Due to the locality of the data involved in this calculation, hardware acceleration may be used to exploit fine-grained parallelism using meshfree libraries such as the Compadre toolkit \cite{compadrev101}.

These quadrature weights may finally be used to evaluate the classical LPS model in Equations \ref{eq:discreteNonlocElasticity} and \ref{eq:discreteNonlocDilitation}:
\begin{equation}\label{eq:discreteNonlocElasticity2}
    \frac{C_\alpha}{m(\delta)} \sum_{j \in B_\delta (\mathbf{x}_i)} \left(\lambda_{ij} - \mu_{ij}\right) K_{ij} \left(\mathbf{x}_j-\mathbf{x}_i\right)\left(\theta_i + \theta_j \right) \omega_{j,i}
\end{equation}
\begin{equation*}
  +  \frac{C_\beta}{m(\delta)}\sum_{j \in B_\delta (\mathbf{x}_i)} \mu_{ij} K_{ij}\frac{\left(\mathbf{x}_j-\mathbf{x}_i\right)\otimes\left(\mathbf{x}_j-\mathbf{x}_i\right)}{\left|\mathbf{x}_j-\mathbf{x}_i\right|^2} \cdot \left(\mathbf{u}_j- \mathbf{u}_i \right) \omega_j = \mathbf{f}_i
\end{equation*}
\begin{equation}\label{eq:discreteNonlocDilitation2}
  \theta(\mathbf{x}) =  \frac{d}{m(\delta)} \sum_{j \in B_\delta (\mathbf{x}_i)}  K_{ij} \left(\mathbf{x}_j-\mathbf{x}_i\right) \cdot  \left(\mathbf{u}_j - \mathbf{u}_i \right) \omega_{j,i}
\end{equation}
and it may be confirmed from direct calculation that the quadrature rule is exact for $\mathbf{u} \in (P_2)^d$ and $\theta \in P_1$.

\section{A modified dilitation definition for damage and boundary effects}\label{ssec:bondCorrection}

As eluded to in Section \ref{ssec:LPS}, the nonlocal definition of dilitation breaks down near domain boundaries and cracks induced by bond-breaking. We briefly recall the bond-breaking procedure here and then introduce a modification to the definition of dilitation to maintain consistency.

In the peridynamic framework, a state is associated with each bond, denoted here as an $(i,j)$ pair, where it is either broken or unbroken. As defined in \cite{silling_2005_2,madenci2014peridynamic}, we modify the quadrature weights obtained in the previous section as 
\begin{equation}
\tilde{\omega}_{j,i} =
  \begin{cases}
    \omega_{j,i}, & \text{ if bond is unbroken }  \\
	0, & \text{ if bond is broken }  .
  \end{cases}
\end{equation}

Typically in peridynamics, a damage model is implemented by prescribing a relationship in which bonds exceding a critical fracture criteria are broken. For example, in the protyptical microelastic brittle (PMB) model which has been used extensively to model brittle materials failing in tension\cite{silling_2005}, bonds exceeding a critical strain are broken, and a strain energy argument is used to calibrate the relationship between nonlocal parameter $\delta$ and a fracture energy.

In this work, we will for simplicity restrict ourselves to statically loaded cracks, whereby a pre-crack is modeled by breaking any bond crossing a fracture surface. In previous work with bond-based peridynamics we have shown that combining OBQ with traditional PMB type models, we were able to preserve asymptotically-compatible convergence to the linear theory with traction-free boundary conditions\cite{trask2019asymptotically}. We will illustrate numerically in Section \ref{ssec:anisoSol} that for both statically loaded cracks and for free surfaces modeled as a fractured surface, we are able to obtain similar results in this state-based case.

\subsection{Numerical scheme}
We modify the classical definition of nonlocal dilitation in Equation \ref{eq:discreteNonlocDilitation2} to enforce that the nonlocal dilitation recovers the local divergence in the case of linear displacement fields, independent of whether the horizon intersects the boundary of the domain. In the spirit of similar techniques employed in the SPH community\cite{oger2007improved}, we accomplish this by modifying Equation \ref{eq:nonlocDilitation} through the introduction of a correction tensor $\mathbf{M}(x)$. 

\begin{equation}\label{eq:continuousNonlocDilitation3}
  \theta^{corr}(\mathbf{x}) =  \frac{d}{m(\delta)} \int_{B_\delta (\mathbf{x})\cap \Omega} K(\left|\mathbf{y}-\mathbf{x}\right|) \left(\mathbf{y}-\mathbf{x}\right) \cdot \mathbf{M}(\mathbf{x})\cdot \left(\mathbf{u}(\mathbf{y}) - \mathbf{u}(\mathbf{x}) \right) d\mathbf{y},
\end{equation}
where, mimicking work by Oger in SPH \cite{oger2007improved}, we define
\begin{equation}\label{eq:continuousdilCorr}
  \mathbf{M}(\mathbf{x}) =  \left[ \frac{d}{m(\delta)} \int_{B_\delta (\mathbf{x})\cap \Omega} K(\left|\mathbf{y}-\mathbf{x}\right|) \left(\mathbf{y}-\mathbf{x}\right) \otimes \left(\mathbf{y}-\mathbf{x}\right)  d\mathbf{y} \right]^{-1}.
\end{equation}
After discretization, we obtain
\begin{equation}\label{eq:discreteNonlocDilitation3}
  \theta^{corr}_i =  \frac{d}{m(\delta)} \sum_{j \in B_\delta (\mathbf{x}_i) \cap \Omega}  K_{ij} \left(\mathbf{x}_j-\mathbf{x}_i\right) \cdot  \mathbf{M}_i^{-1} \cdot \left(\mathbf{u}_j - \mathbf{u}_i \right) \omega_{j,i},
\end{equation}
and
\begin{equation}\label{eq:dilCorr}
\mathbf{M}_i = \left[\frac{d}{m(\delta)} \sum_{j \in B_\delta (\mathbf{x}_i) \cap \Omega}  K_{ij} \left(\mathbf{x}_j-\mathbf{x}_i\right) \otimes  \left(\mathbf{x}_j - \mathbf{x}_i \right) \omega_{j,i}\right].
\end{equation}
Evaluation of the correction tensor requires that Equation \ref{eq:dilCorr} be invertible; this is true as long as there are at least $d$ non-collinear bonds within the horizon. In some settings, such as violent dynamic fracture, for a given particle all bonds may break, leaving an isolated particle. In this case the matrix inverse $\mathbf{M}_i^{-1}$ may be replaced with a pseudo-inverse $\mathbf{M}_i^+$ to improve robustness of the scheme.

\section{Numerical results}\label{sec:numerics}
For all solutions presented, we will take $\Omega$ to be the unit square. Dirichlet conditions are imposed on a layer of size $\delta$ around the boundary of $\Omega$. We generate a lattice of spacing $h$ within the domain, and then perturb each point by a uniformly distributed random variable of magnitude $0.2\, h$ to ensure a quasi-uniform pointset with no symmetry. For all problems we use a quadratic reproducing space in the quadrature rule, to match the second order convergence of the nonlocal model to classical elasticity, and scale the kernel radius as $\delta = 3.5 h$ to ensure sufficient neighbors to obtain polynomial unisolvency. 

To model the traction free interface in problems where $\partial \Omega_N \neq \emptyset$, we will break any bond intersecting the Neumann boundary. We highlight that quadrature weights are computed on the perturbed lattice of points before bonds are broken. To visualize broken bonds, we will post-process the following damage measure.
$$d_i = 1 - \sum_{j \in B_\delta(\mathbf{x}_i)} \frac{\tilde{\omega}_{j,i}}{{\omega}_{j,i}}$$

Error is computed with respect to the root-mean-square norm
$$||f||_{\ell_2} = \sqrt{\frac{\sum_i f_i^2}{N}},$$
which for quasi-uniform pointsets may be shown to be equivalent to the full $L_2(\Omega)$ norm \cite{wendland2004scattered}, i.e. there exists $C>0$ such that
$$||f||_{L_2(\Omega)} := \int_\Omega |f|^2 dx \leq C ||f||_{\ell_2}.$$

For the governing equations, we will obtain a block structured matrix with respect to the displacements and dilitation. We solve this system with a Schur-complement preconditioner described in a previous work\cite{trask2016compact}. While an in-depth study of fast solvers for these problems is beyond the scope of this study, we note that the above strategy yields a roughly $O(N log N)$ scaling for problems with moderate Poisson ratio (i.e. far from the incompressible limit).

\subsection{Manufactured solutions}\label{ssec:manSol}
We consider the following two manufactured solutions with corresponding forcing term.

\begin{align}\label{eq:man1}
\mathbf{u}_1 = \left<x^2,4 y^2\right>\\
\mathbf{f}_1 = \left<2 \lambda + 4 \mu, 8 \lambda + 16 \mu\right>,
\end{align}
and
\begin{align}\label{eq:man2}
\mathbf{u}_2 = \left<\sin x \sin y, -\cos x \cos y\right>\\
\mathbf{f}_2 = -2 \left(\lambda + 2 \mu\right) \left<\sin x \sin y, -\cos x \cos y\right>.
\end{align}
The first provides a patch test illustrating the ability of the discretization to exactly reproduce polynomial solutions, while the second illustrates the optimal second order convergence\cite{trask2019asymptotically}. We summarize results of a refinement study in Table \ref{tab:manu}, considering solutions to both Equations \ref{eq:man1} and \ref{eq:man2} for $\mu = \lambda = \frac12$ and solutions to Equation \ref{eq:man2} for $\mu$ and $\lambda$ corresponding to a near-incompressible material with Poisson ratio $\nu = 0.495$ and Youngs modulus $E=1$.

Note that for this problem the dilitation correction tensor (Equation \ref{eq:dilCorr}) is the identity; because each horizon is completely supported within the domain, the dilitation correction does not do anything, and the dilitation reverts to the classical definition.

\begin{table}[]\label{tab:manu}
\center
\begin{tabular}{cccc}
\textbf{N} & $||u_1-u_1^{ex}||_{\ell_2}$ & $||u_2-u_2^{ex}||_{\ell_2}$ & $||u_{2,inc}-u_{2,inc}^{ex}||_{\ell_2}$ \\ \hline
$24^2$     & 4.11e-14                    & 0.02207  & 0.13057                    \\
$48^2$     & 2.47e-13                    & 0.00506  & 0.02597                \\
$96^2$     & 2.69e-12                    & 0.00117  & 0.00632                    \\
$192^2$    & 4.01e-13                    & 0.00028  & 0.00158                   
\end{tabular}
\caption{ Convergence study for manufactured solution problem. $u_1$ corresponds to quadratic patch test, $u_2$ to smooth solution, and $u_{2,inc}$ to smooth solution near incompressible limit.}
\end{table}

\subsection{Isotropic materials with free surfaces and fracture}\label{ssec:isoSol}

We consider now convergence to the analytic solution of a hole of radius $a$ in an infinite plate subject to far field uniaxial tension $T$. The solution to this problem is classical and may be obtained via a Michell solution. Under a plane stress assumption, the displacement field is given by

\begin{align}\label{eq:holeSol}
u_x(r,\theta) = \frac{T a}{8 \mu} \left[ \frac{r}{a} (\kappa + 1) \cos \theta + \frac{2 a}{r}\left((1+\kappa) \cos \theta+\cos 3\theta \right) - \frac{2 a^3}{r^3} \cos 3 \theta \right]\\
u_y(r,\theta) = \frac{T a}{8 \mu} \left[ \frac{r}{a} (\kappa - 3) \sin \theta + \frac{2 a}{r}\left((1-\kappa) \sin \theta+\sin 3\theta \right) - \frac{2 a^3}{r^3} \sin 3 \theta \right],
\end{align}
where the coefficient
\begin{equation}
\kappa = 
\begin{cases}
3 - 4 \nu           \quad \text{ for plane strain}\\
\frac{3-\nu}{1+\nu} \quad \text{ for plane stress}
\end{cases}
\end{equation}

This case provides an important benchmark because it allows verification of asymptotic compatibility in the presence of a curvilinear free surface induced by bond breaking; historically surface effects in this scenario have proven problematic for state-based peridynamics. In this scenario, the dilitation correction (Equation \ref{eq:dilCorr}) does take a non-trivial value.

\begin{figure}[t!]
   \centering
   \includegraphics[width=0.47\textwidth]{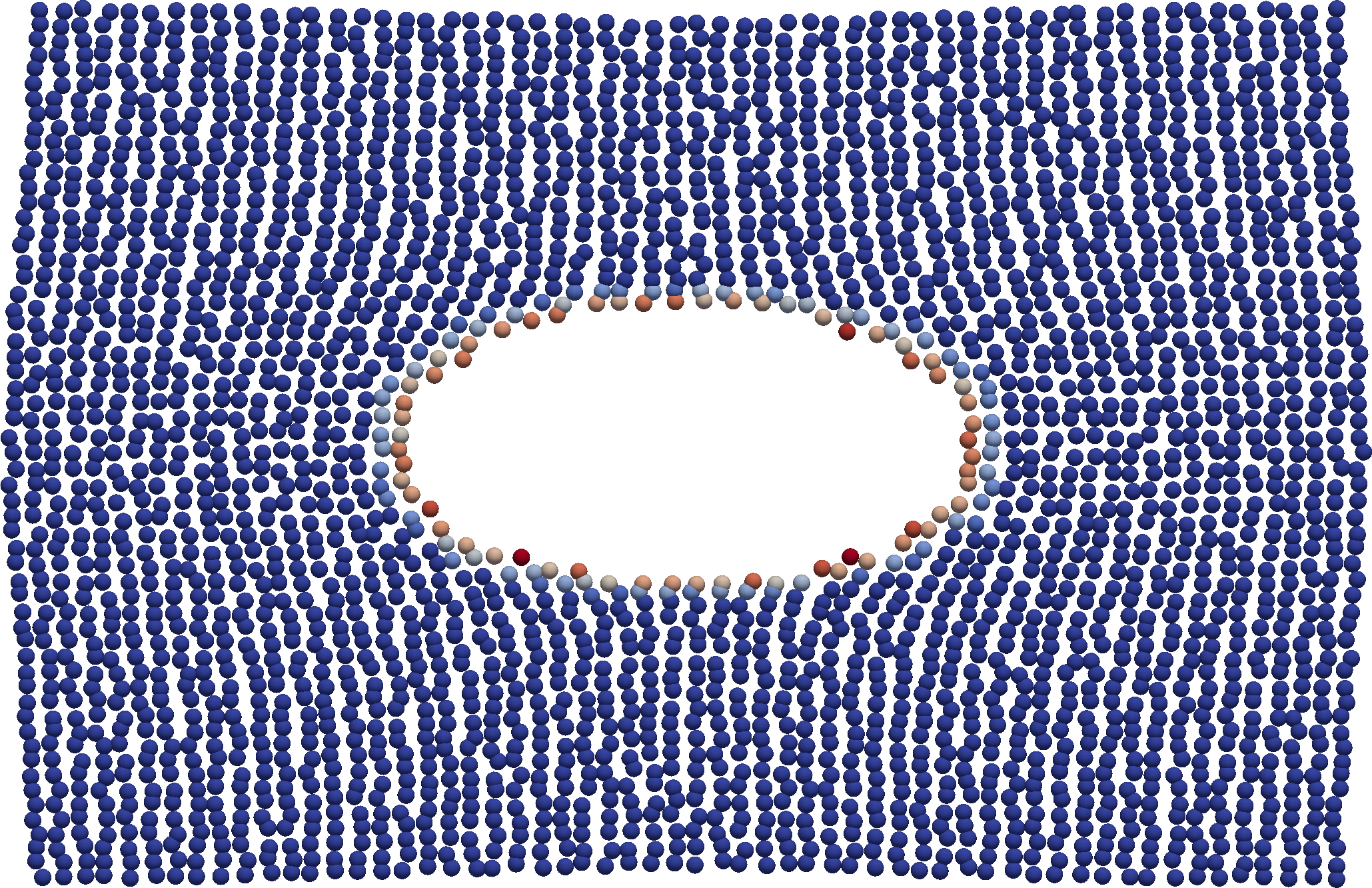}
   \includegraphics[width=0.47\textwidth]{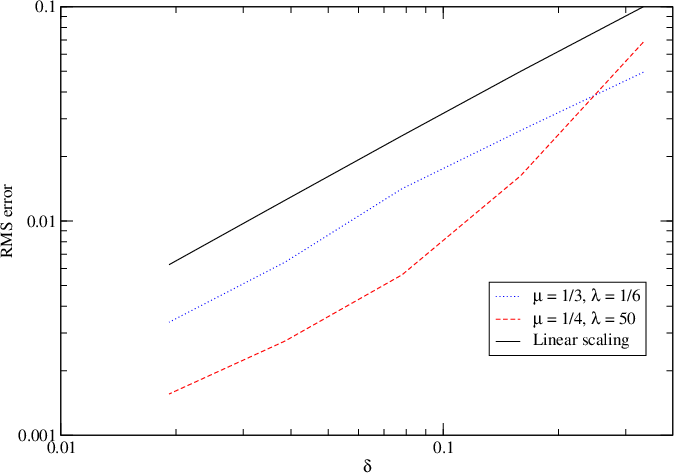}\\
   \caption{\textbf{Left:} Damage in deformed configuration, illustrating free surface treatment by breaking bonds. \textbf{Right:} First-order convergence of displacements in RMS norm for $\lambda$ and $\mu$ chosen both close to and far from the incompressible limit.}
   \label{fig:holeInPlate}
 \end{figure}

To set the problem up, we impose the analytic solution as a Dirichlet condition on the nonlocal collar around the perimeter of the unit square. We then break all bonds crossing the circle of radius $a = 0.2$. In Figure \ref{fig:holeInPlate}, we visualize the free surface by plotting damage in a deformed configuration, and provide a convergence plot demonstrating asymptotic compatibility both far from and close to the incompressible limit.

\subsection{Composite materials with internal interfaces}\label{ssec:anisoSol}

We consider now a solution to a hydrostatically loaded cylindrical inclusion of radius $a$ in an infinite plate. We denote the interior of the inclusion as $\Omega_1$ and the exterior as $\Omega_2$, with associated constant material properties $\left(\mu_1,\lambda_1\right)$ and $\left(\mu_2,\lambda_2\right)$. 

Assuming a far-field hydrostatic stress $P_\infty$ and plane strain conditions, we define the coefficients

\begin{align}
C_A &=   \frac{P_\infty}{2(\lambda_1 + \mu_1)}\\
C_B &=   \frac{P_\infty\left(\lambda_1 + \mu_1 +\mu_2\right)}{2(\lambda_1 + \mu_1) (\lambda_2 + 2 \mu_2)}\\
C_C &= - \frac{P_\infty a^2\left(\lambda_1-\lambda_2 + \mu_1 -\mu_2\right)}{2(\lambda_1 + \mu_1) (\lambda_2 + 2 \mu_2)}
\end{align}

and the analytic solution for the displacement field in spherical coordinates is given by 

\begin{align}
u_r &=
\begin{cases}
C_a r         \quad \mathbf{x} \in \Omega_1\\
C_b r + C_c/r \quad \mathbf{x} \in \Omega_2
\end{cases}\\
u_\theta &= 0.
\end{align}

We use this solution to assess the stability of the method in the vicinity of large jumps in material properties. Note that the consistency conditions derived only guarantee asymptotic compatibility under the assumption of an isotropic material; this benchmark thus explores the applicability of the approach beyond the guarantees of the approximation theory in \cite{trask2019asymptotically}.

We first investigate whether the discretization is asymptotically compatible. To set this up, we take again $a = 0.2$, impose a jump in the bulk modulus ($K_1 = 2$, $K_1 = 1$), and consider again two scenarios corresponding to far- and near- incompressible materials: $\nu_1 = \nu_2 = 0.25$, and $\nu_1 = 0.49$ and $\nu_2 = 0.25$. We consider both a uniform Cartesian grid and the perturbed quasi-uniform grid discussed previously. The results are presented in Table \ref{tab:inclusion}. It is apparent that for a uniform grid, for both Poisson ratios we obtain asymptotically compatible convergence to the local solution. However, for the non-uniform grid, the convergence is seen to stagnate in the compressible material after obtaining four digits of accuracy. While this is sufficient accuracy for many engineering problems, this suggests that a Cartesian grid is necessary if asymptotic compatibility is to be ensured.

\begin{figure}[t!]   \label{fig:inclusion}
   \centering
   \includegraphics[width=0.32\textwidth]{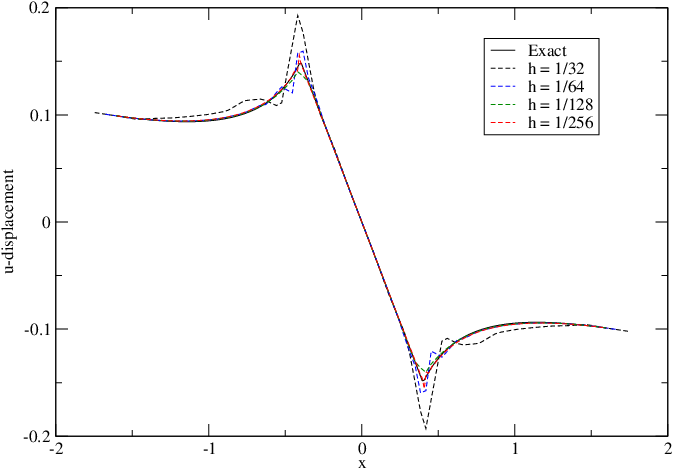}
   \includegraphics[width=0.32\textwidth]{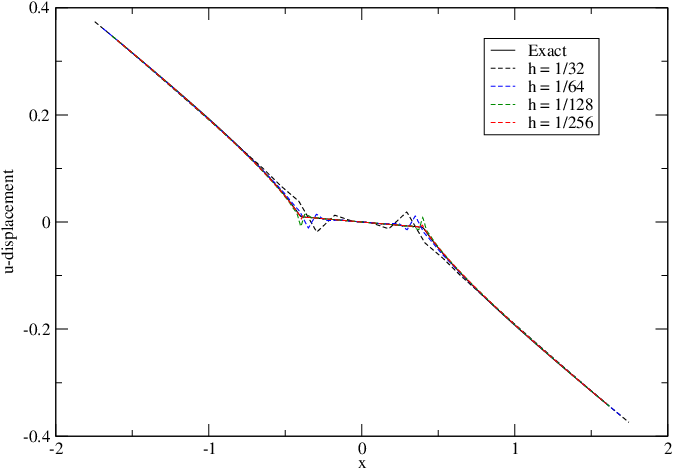}
   \includegraphics[width=0.32\textwidth]{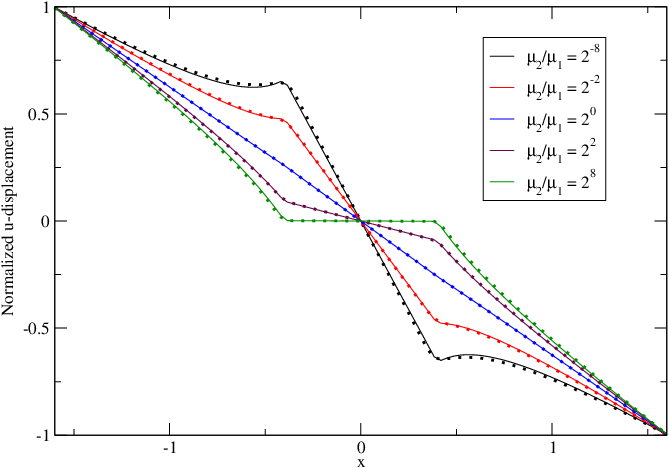}
   \caption{\textbf{Left:} Convergence to analytic solution for a soft inclusion ($\mu_2/\mu_1 = 1/64$). \textbf{Center:} Convergence to analytic solution for a stiff inclusion ($\mu_2/\mu_1 = 64$). \textbf{Right:} For a fixed resolution of $64^2$ points, reproduction of analytic solution for a wide range of $\mu_2/\mu_1 \in \left(2^{-8},2^8\right)$. Solid line corresponds to analytic solution, while dots correspond to numerical result.}
 \end{figure}

We next investigate the stability of the approach over a large range of material parameters. To do this, we fix the Poisson ratio in both phases to $\nu_1 = \nu_2 = \frac14$ and impose a jump in the shear modulus of $\mu_2/\mu_1 = K$, for $K \in \left\{2^{-8}, 2^8\right\}$. In Figure \ref{fig:inclusion}, we plot a profile of the x-component of displacement along the $y = 0$ line. We demonstrate convergence for both a stiff inclusion ($\mu_1 = 64 \mu_2$), a soft inclusion ($\mu_2 = 64 \mu_1$), and then illustrate that we reproduce well the displacement for a wide range of parameters. Convergence results are documented in Table \ref{tab:inclusion}.

\begin{table}[]
\center
\begin{tabular}{ccccc}
\multicolumn{1}{l}{} & \multicolumn{2}{c}{Uniform}     & \multicolumn{2}{c}{Nonuniform}                                          \\
\textbf{N}           & $\nu_2 = 0.25$ & $\nu_2 = 0.49$ & \multicolumn{1}{l}{$\nu_2 = 0.25$} & \multicolumn{1}{l}{$\nu_2 = 0.49$} \\ \cline{1-5}
$16^2$               & 0.00569        & 0.04463        & 0.00661                            & 0.04629                            \\
$32^2$               & 0.00201        & 0.03926        & 0.00242                            & 0.03941                            \\
$64^2$               & 0.00099        & 0.02260        & 0.00144                            & 0.02304                            \\
$128^2$              & 0.00045        & 0.01205        & 0.00055                            & 0.01211                            \\
$256^2$              & 0.00023        & 0.00595        & 0.00044                            & 0.00614                           
\end{tabular}
\caption{Convergence in RMS norm for inclusion problem considering uniform and perturbed particle distributions both close-to and far-from the near incompressible limit.}
\label{tab:inclusion}
\end{table}

\section{Conclusions}
We have presented a modification to the definition of dilitation in state-based peridynamics, so that the classical linearly elastic peridynamic solid recovers the AC limit to classical elasticity as the nonlocal and discretization lengthscales are reduced at the same rate near boundaries. By construction, the approach is able to pass a patch test to machine precision. We show that for isotropic materials, this convergence holds in the presence of curvilinear free surfaces, and thus provides a scheme in which the bond-breaking techniques popular in the peridynamics community may be consistently applied to linearly elastic materials. We have also shown that for composite materials the approach is only able to preserve the AC limit for Cartesian particle arrangements, although reasonable accuracy is achieved for quasi-uniform pointsets.

The ubiquity of examples where peridynamics is used as a modelling platform to describe physical processes involving microstructure and fracture is a testimony to the flexibility and descriptiveness of the theory; historically however these works often fail to demonstrate convergence quantitatively. We posit that discretizations with quantitatively benchmarked notions of consistency will lend credibility to nonlocal models overall and are fundamental to promoting the broader adoption of nonlocal methods. In light of this, while many works have pursued consistent discretization of state-based peridynamics in the past, to our knowledge this work marks the first AC meshfree discretization of nonlocal elasticity with bond-breaking. In future work, we will pursue the application of this scheme to study nonlinear elastoplastic mechanics and ductile fracture.

\section{Acknowledgements}

    This work was supported by the (other funding source references), and the Laboratory Directed Research and Development program at Sandia National Laboratories. Sandia National Laboratories is a multi-program laboratory managed and operated by National Technology and Engineering Solutions of Sandia, LLC., a wholly owned subsidiary of Honeywell International, Inc., for the U.S. Department of Energy's National Nuclear Security Administration under contract DE-NA-0003525. 

\bibliographystyle{unsrt}
\bibliography{statebased}

\end{document}